\documentclass{notices}
\usepackage{graphicx, amsmath, textcomp, amsrefs}
\title{
Joseph Carter Corbin: Arkansas's ``Profound Mathematician"
}

\author{
  Jesse Leo Kass
  \affil{
    Jesse Leo Kass is an assistant professor at UC Santa Cruz. His email address is jkass@umich.edu.
    }
}

\begin{document}

\maketitle

While Joseph Carter Corbin  (1833--1911) is a celebrated figure remembered for his impact in higher education and politics within the state of Arkansas,\footnote{Corbin's life is the subject of book \cite{corbinBioBook} and an article \cite{rothrock}. The present account of Corbin's life draws heavily from those publications, and we provide citations only for those facts not documented in those sources.} he is absent from many modern historical accounts of American mathematicians.\footnote{Two notable exceptions which discuss Corbin's publications in the {\it American Mathematical Monthly} are \cites{means59, MAA_article15}.}   However, during his lifetime Corbin was very much part of the mathematical community. He was described as a ``Profound Mathematician" in {\it Men of Mark}, an influential anthology of accomplished Black men \cite{menOfMark}*{pp. 829--832}. He regularly contributed to mathematical periodicals like the {\it American Mathematical Monthly}.  Upon his death,  the  {\it Monthly} commemorated him by  publishing his obituary \cite{corbinObit}. 

Corbin's omission from historical accounts of American mathematicians is neither surprising nor unusual. Most  mathematicians of his generation are omitted because accounts tend to focus on  research activities. Only towards the end of Corbin's life did significant numbers of American math professors begin to do research.

Certainly Corbin is worthy of inclusion in the history of American mathematics. His life was remarkable. The son of freed slaves, Corbin was the founding head of the University of Arkansas at Pine Bluff, Arkansas's public Historically Black University. While working there,  for over three decades, he regularly published in mathematical periodicals. 

With the goal of introducing Corbin to modern mathematicians, we survey his life in this article. We focus on his experiences in higher education and mathematics and, in particular, give an overview of his mathematical publications.

\subsection*{Early Life}
Corbin was born in  Ohio in 1833. His parents had been enslaved in Virginia but had moved to Ohio about a decade before Corbin's birth. Little is known about their move, but it was a common one: most Black Ohioans\footnote{We will use the adjective ``Black" to refer to enslaved Americans and their descendants. This is somewhat ahistorical. Corbin was often regarded as being of mixed race even though both parents had been enslaved. For example,  federal census takers recorded Corbin's race as  ``mulatto" in every census between 1850 and 1910 except for 1900 (when it was recorded as ``Black") and possibly 1890 (this record is unavailable). See also the discussion of Corbin's time at Ohio University in the section titled ``Education."}  had come from the bordering slave states of Virginia and  Kentucky. Corbin's mother, Susan, had moved after being emancipated by her enslaver. It is unknown whether the  father, William, was also emancipated or if he had fled enslavement. 

Even though Corbin never personally experienced enslavement, his life was significantly constrained by state laws and social practice. For example, Black Ohioans were not allowed to allowed to attend public schools, and many private schools were racially segregated. Despite the obstacles facing the family, Corbin's parents were able to achieve considerable personal success;  both were literate,  and William worked as a ``sewer of clothes" (a skilled profession).\footnote{See William's entry in the 1850 U.S.~census. He was described as a ``porter" in the 1860 census.}

\subsection*{Education}

When Corbin was growing up, his family lived in Chillicothe, a mid-sized town of a few thousand people in southern Ohio which was a regional center for Black life (Blacks residents made up 10\% of Chillicothe's population but only 1\%  of the state's). Corbin attended school in Chillicothe until he was fifteen years old (around 1848). He then moved to Louisville, Kentucky to take advantage of the greater educational opportunities afforded by a larger city. He returned to Ohio in 1850 to attend Ohio University in Athens.

Black students had studied at OU since the 1820s, and Corbin was the third student  to enroll. While the presence of  Black students on college campuses was sometimes a source of tension,  Corbin did not seem to have attracted any special attention. It is unclear whether Corbin was even regarded as Black by faculty and students. For example, an 1853 newspaper article on the annual commencement exercises at Ohio University mentions Corbin by name but makes no reference to his race, a standard journalistic practice at this time \cite{newspaper53}. Similarly, an 1885 list of university  alumni \cite{catalogue18885}*{p.~86} describes  Corbin simply  as a teacher in Louisville\footnote{The description of Corbin as a teacher in Louisville suggests that his entry in the list was based on old  information collected shortly after Corbin graduated (during the 1850s or 1860s, when he was in Louisville). At the time the list was published, Corbin had been living in Arkansas for over a decade.} but describes another alumnus  (John Newton Templeton) as ``the only Alumnus .~.~.~of African descent" \cite{catalogue18885}*{p.~78}.

When Corbin arrived at OU, the university was functioning as a standard 19th century U.S.~university. It offered a B.A.~degree for students who completed a  4-year sequence of college courses. The university also maintained a Preparatory Department which offered a 2-year sequence designed to prepare students for the college courses.\footnote{While uncommon today, many 19th century universities maintained similar departments. Indeed, doing so was often necessary as schooling along the lines of modern K–12 education was often only available on a limited basis.} At the start of Corbin's first semester, enrollment stood at  64 students,  divided more or less evenly between college and preparatory students. They were taught by five professors, each being responsible for all coursework in a given subject. The graduation requirements differed from those of a modern American university in that students did not select majors. Instead, all college students took the same fixed sequence of courses. 

The required coursework heavily emphasized the study of Greek and Roman literature, but it  also included a mathematics curriculum that covered  algebra and elementary geometry in the first year; trigonometry, analytic geometry, and differential calculus in the second year; and integral calculus in the third year.\footnote{Corbin's analytic geometry and calculus classes were taught using  Albert E.~Church's books  {\it Elements of the Differential and Integral Calculus} and {\it Elements of Analytical Geometry}, while the algebra and geometry classes were taught using textbooks by Charles Davies \cite{OUCatalogue51, OUCatalogue53}. Algebra was taught from {\it Elements of Algebra: Translated from the French of M. Bourdon}. It is not entirely clear which geometry textbook were used as Davies wrote several. The assigned text is listed as ``Davies' Legendre" in the 1851-50 university catalogue and as ``Davies" in 1852-53 university catalogue. The first reference is to {\it Elements of Geometry and Trigonometry from the works of A. M. Legendre}, but the second reference could either be to that book or to another one such as {\it Elements of Descriptive Geometry}.}

	  At OU, Corbin was taught by two different math professors,  William J.~Hoge and Addison Ballard.  As was common for American math professors at the time, neither  Hoge nor Ballard had any specialized training in math.\footnote{See \cite{butonBook}*{pp.~659-669} for a more detailed discussion of the education of a typical math professor in nineteenth century America. Formal advanced education in mathematics was out of reach for most  as it typically necessitated studying abroad. American universities only began awarding earned doctorate degrees in 1861.}  Corbin's first math classes were taught by Hoge, an alumnus of the university who had graduated  with an undergraduate degree in 1843 \cite[p.~84]{catalogue18885}. Hoge left the university while Corbin was  still a student (in 1851), and he was replaced  by  Ballard. Ballard had received a bachelor's degree from Williams College and had been serving as OU's  Latin Professor since 1848 \cite[pp.~213-214]{williamsBook}.

\subsection*{Career}
When Corbin graduated from Ohio University in 1853, his formal education came to an end,\footnote{Corbin would later receive two A.M.~degrees from OU and a Ph.D.~from an unknown Baptist college, but these were honorary degrees.}  and he spent the next two decades working as a teacher, bank clerk, and newspaper editor in  Louisville and Cincinnati.

In 1872, Corbin moved to the former slave state of Arkansas. This was not an uncommon move for talented Black men like Corbin, as the political changes brought about by the Civil War created unprecedented opportunities for change. After  working briefly  as a newspaper reporter and a post office clerk, Corbin was elected to oversee Arkansas's public education system as superintendent of public instruction.\footnote{Corbin defeated the incumbent, Thomas Smith, a white Union army veteran from Pennsylvania. Both Smith and Corbin were Republicans, but they were aligned with rival groups within the party \cite{Kennan61}.} His election made Corbin the head of the  board of trustees for the state's newly founded public university --- Arkansas Industrial University (now the University of Arkansas).\footnote{Prior to the Civil War, Arkansas did not maintain a public university, although there were a few  private universities. During the war, U.S. Congress encouraged state governments to create public universities by passing the Morrill Land-Grants Acts, which provided states with funding for public universities through land sales. Arkansas Industrial University was created to fulfill the terms of the act so that the state could receive funding.}

The  university was technically ``open to all, without regard to race, sex, or sect,"\footnote{The quotation is taken from a statement by the chair of the executive committee of the board of trustees in an April 1873 report. Corbin probably did not play a major role in the decision to open the university to Black students. He was not on the executive committee, and the committee was charged with writing the report before he joined the board.}  but Black students were not allowed to attend classes with white students. Instead, they received private tutoring.\footnote{It is unclear exactly how many Black students attended, although the numbers were certainly small. A 1910 history of the university states that only one Black student applied for admission \cite[p.~97]{ARhistoryBook}, but a 1922 letter by the interim president's wife states that two or three applied \cite[p.~3]{nichols}. At least one Black student, James McGahee, attended the university from 1872-73.} To address the issue of college education for Black Arkansans, the board of trustees, including Corbin, petitioned the state legislature to create a school for educating Black teachers. The legislature responded favorably to the petition and created the Branch Normal College (now the University of Arkansas at Pine Bluff).\footnote{The legislative act  made no explicit reference to race, but it was tacitly understood that  only Black students would be allowed to attend the college. For example, the act directed the board to select a site for the college in the southeastern part of the state (where most of Arkansas's Black population lived) that was convenient for the ``poorer classes" (understood to be a euphemism for Blacks Arkansans).} 

  Corbin's  term in office was cut short when Conservative Democrats, who were largely hostile towards Black Arkansans,  gained control over the state government and vacated many political offices, including Corbin's. Corbin left the state and moved to Jefferson City, Missouri, where he taught at the newly founded Lincoln Institute (now University).\footnote{At the time, Lincoln was a small school of about 50 students, mostly the children of freedpersons \cite{newspaper74, newspaper75a, newspaper75b}.} However, after a brief time there, he returned to Arkansas  to serve as the first principal of the Branch Normal College.

 Since the 1873 legislative act, no real progress had been made in establishing the college. Under state law, the Branch Normal College was supposed to provide the same education as the Normal Department at Arkansas Industrial University.  A normal department or school was a common 19th-century institution that  offered  postsecondary education for  school teachers which typically involved  a truncated college curriculum along with specialized coursework in pedagogy.  However, achieving parity with Arkansas Industrial University was an unrealistic goal for the Branch Normal College as many of the college's students were facing rural poverty and political disenfranchisement. To further complicate the situation, Corbin was given few resources --- the college had no permanent facilities and no faculty aside from Corbin, who was even responsible for menial tasks like cleaning classrooms. Nevertheless, Corbin was able to make major improvements. By the 1880s, a college building and dormitories had been constructed, and the faculty had increased to three.  By 1883, student enrollment had reached over 200 students.

The Branch Normal College offered an academic education through three different departments: (a) the Preparatory Department which offered remedial classes to prepare students for college work,  (b) the Normal Department which awarded a Licentiate of Instruction (a teaching certificate), and (c) the Collegiate Department which awarded a Bachelors of Arts. The majority of students were enrolled in the Preparatory Department. The preparatory curriculum  consisted of three years of courses on subjects such as grammar, geography, drawing, and math. The math education consisted of six terms of elementary arithmetic and one of algebra.\footnote{According to the 1882 college catalogue \cite{UofArkansasCatalogue82}*{p.~91-92}, arithmetic was taught  using books from  Robinson's Series of Mathematics, a popular series of textbooks. It is unclear precisely which books were used as they are only described as ``Arithmetic, Robinson's Shorter Course" and ``Arithmetic, Robinson's." Algebra was taught from William G.~Peck's {\it Manual of Algebra}.}

Students in the Normal Department completed two years of remedial work followed by two years of college courses. In addition to classes on Latin, English grammar, and pedagogy, students studied algebra and geometry in their third year and then plane trigonometry in their fourth.\footnote{The 1882 college catalogue \cite{UofArkansasCatalogue82}*{p.~91-92}   states that geometry  was taught from {\it Wentworth's Solid Geometry} by George Wentworth, and trigonometry from Edward Olney's {\it Elements of Trigonometry, Plane and Spherical}. The algebra textbook is described as ``Thompson's Algebra,"  but the present author has been unable to identify this book.} Students in the Collegiate Department completed the courses taken by the normal students and then two additional years of college classes in subjects like foreign languages, English literature, and mathematics. The math classes offered in the last two years were a class on analytic geometry and an optional course on calculus.\footnote{The description of courses is taken from the 1893 college catalogue  \cite{UofArkansasCatalogue93}. The textbooks used are not listed.}

The curriculum that Corbin promoted at the Branch Normal College came into conflict with an increasingly influential educational philosophy that focused on industrial education. This philosophy endorsed only a very limited academic education for Black students and instead emphasized the value of manual labor as a means to impart values like self-discipline and self-reliance. Many of those who promoted industrial education scorned  Black students who saw academics as a way of obtaining personal advancement or fulfillment. Booker T.~Washington, the most prominent promotor of industrial education, expressed the views of many in his autobiography. In an account of teaching in Alabama, he wrote that,  ``The students who came [to his school] first seemed to be fond of memorizing long and complicated `rules' in grammar and mathematics, but had little thought or knowledge of applying these rules to the everyday affairs of their life" \cite{bookerBook}*{p.~122}. To illustrate his point, he offered the example of a student whose fondness for math led him to study methods for computing cubic roots. Washington set students like him to activities that he deemed more appropriate:  farm work and learning how to properly set a dinner-table.

Seeking to implement industrial education in Arkansas, state legislators in the 1890s created a Department of Mechanical Arts at the Branch Normal and hired as department head William S.~Harris, a white man from Virginia. He was effectively made head of the college by the board of trustees when they transferred key responsibilities\footnote{For example, Harris was made responsible for admitting students and collecting fees.}  from Corbin to Harris.

Corbin nominally remained the principal of the Branch Normal for approximately another decade. In 1902, the board of trustees replaced him by hiring Isaac Fisher, a recent graduate of Booker T.~Washington's Tuskegee Institute.  Fisher's hire was part of a plan to expand industrial education at the Branch Normal, but Fisher's efforts to the transform the college faltered in the face of mixed support from trustees and strong opposition by the Black community. While Corbin's removal was a major set-back to efforts to run the Branch Normal  as an academic institution, the college continued to offer advanced courses taught by Corbin's former student James C.~Smith.

After his removal from the college, Corbin remained in Pine Bluff\footnote{Corbin's obituary in \cite{corbinObit} states that, from 1902 to January 1904, he served as president of Ouachita Baptist College (now University) in Camden. While newspaper records from the time record Corbin's presence in Camden, the rest of this statement seems to be erroneous. During this period, the college's president was one John W.~Conger, not Corbin. The college was also not located in Camden. It has always been located in Arkadelphia, a town about 40 miles away. In fact, the college has no record of an association with Corbin. His employment, in any capacity, would have been highly unusual as  Ouachita was a whites-only college during the 19th century. The author thanks Ouachita Professor and University Archivist Lisa Speer for correspondence on this issue.} and served as principal of Merrill Public School (the local Black high school) until his death in 1911, at the age of seventy-seven.

\subsection*{Corbin's mathematical work}
 Corbin first published in a mathematical periodical in 1882, when he was almost fifty years old. At this point in his life, Corbin had achieved remarkable success in politics and education, but this publication is the first record of his intellectual engagement with mathematics. This was a point of pride for Corbin, but it is unclear if these publications were the expression of a life-long interest or one developed late in life. His publication record shows  both a level of sophistication and a depth of knowledge that extended beyond both what he'd been taught at university and what he was teaching at the Branch Normal College. For example, the calculus textbook he used as a student contained essentially no proofs and even omitted basic definitions like the definition of the definite integral as a limit of Riemann sums. Only the geometry courses he took went beyond mechanical work and included detailed proofs as part of a treatment of plane geometry.

The periodicals that Corbin contributed to were similar to, and in fact included, the {\it American Mathematical Monthly} which is currently published by the Mathematical Association of America.\footnote{The MAA did not exist during Corbin's lifetime, and the {\it Monthly} was privately published by the math professor Benjamin Finkel.} A typical periodical included a problems and solutions section together with book reviews and expository articles. Compared to today's research journals, these periodicals served a very different audience and performed a different function. They were not intended as long-term records of original mathematical research. Instead, they were published to stimulate activity among people who otherwise had few outlets for their mathematical interests. Most readers were individuals studying math in relative isolation, with only limited access to mathematical literature and other resources. They included both amateur mathematicians and math professors, who were often the sole mathematicians at their universities.\footnote{For example, in addition to university professors, the contributors to the first volume of the {\it Monthly} included a counsellor at law (John Doman) and an examiner for the U.S.~Civil Service Commission (Theodore L.~DeLand).} In an era when mathematicians were only beginning to organize themselves into professional societies, these periodicals performed a similar function in fostering a mathematical community.\footnote{While the MAA did not exist during Corbin's lifetime, the  American Mathematical Society was formed  six years after Corbin's first publication, in 1888. It was originally a New York-based society named the New York Mathematical Society, but it was given its current name and reorganized as a national society in 1894. Corbin was never an AMS member, and in general, the AMS was run as an elite and exclusive organization during his lifetime. Admission to the society required receiving a majority of ballots by members in an election after being proposed by two members and recommended by the AMS Council. In contrast, the periodicals Corbin published in were run in a democratic spirit. For example, in the introduction to the inaugural issue of the {\it Mathematical Visitor}, the editor wrote that he aimed to produce an ``unpretending periodical" that inspires a love of mathematics, and he invited contributions from ``professors, teachers, students, and all lovers of the `bewitching science'" \cite{mathVisitor78}' }

Corbin's first published problem  is the following one:
	
\begin{quote}
	Separate the fraction $a/b$ into two parts,  $c/d$ and $e/f$ such that $c+e=d+f$.
\end{quote}
This problem  appeared in a 1882 issue of {\it The Mathematical Magazine}  \cite{corbin82}. The subsequent issue contained two methods for finding the solution $c/d=(a-b-1)/(b-1)$ and $e/f = (b^2+b-a)/(b^2-b)$.\footnote{Although this is the only solution that was published, there are others. For example, $c/d=(a n - b n - 1)/(n b - 1)$ and $e/d=(n b^2 + b - a)/( n b^2 - b)$ is a solution for any integer $n \in \mathbf{Z}$.} Close consideration of the solutions suggests that caution needs to be taken when assigning authorship. Each solution was attributed to two different men, but these men (four in total) probably did not write the published text describing the methods.  Consider the first method.\footnote{The method is to solve for $e/f$ as $a/b- c/d = (a d-b c)/(a d)$, substitute this expression into the equation $c+e=d+f$, and then solve for $c/d$ \cite{corbin82b}*{p.~65}.} Although it is attributed  to P.~M.~McKay and  William E.~Heal, it is unlikely that they jointly wrote the text (one lived in Arkansas, the other in Indiana; there is also no other record indicating collaboration). It is more likely the case that the editor simply summarized two submissions that he received.  

Even when authorship of a work is clear, it did not always represent original ideas. It was not uncommon to submit routine problems  taken from textbooks. In the context of the intellectual environment of the 1890s, this is unsurprising. Textbook problems were interesting to readers as many had little or no access to the books themselves.  For example, the following problem, which Corbin published in a 1898 issue of the {\it Monthly} \cite{corbin98a}, appears to have been taken from the textbook {\it Treatise of Ordinary and Partial Differential Equations}\footnote{Compare the text of Corbin's problem with Example IX.7 on page 104 in the chapter ``Linear Equations: Constant Coefficients" of Johnson's book. This first observation was made by Walter Hugh Drane in his published solution \cite{corbin98b}.} by William Woolsey Johnson: 
\begin{quote}		
	Form the differential equation of the third order, of which $y = c_1 e^{2 x} + c_2 e^{-3 x} + c_{3} e^{x}$ is the complete primitive.
\end{quote}
The differential equation is $y^{(3)} - 7 y'+6y=0$, and it can be solved using the techniques from Johnson's book --- standard techniques still taught to undergraduate students. Despite its routine nature, this problem appears to have interested readers (a number of solutions were submitted \cites{corbin98b, corbin98c}). Among those whose solutions were published was Princeton University professor Edgar Odell Lovett, who gave a complete solution even though he remarked that it was ``a familiar one to students of differential equations." 

Lovett's solution is particularly notable as it demonstrates that, despite their elementary nature, Corbin's publications were even reaching American research mathematicians with an international reputation. Lovett  had completed a doctorate in Germany, published in research journals, and was a member of several European professional societies.

Corbin's differential equation problem must have been the product of significant self-study as the subject was not taught at Ohio University when Corbin was a student. Other problems drew on the plane geometry that Corbin had been taught at Ohio University. One such example is the following  problem \cite{corbin07} which Corbin published in a 1907 issue of the {\it Monthly}:
\begin{quote}
	In triangle $ABC$, the triangle $DEF$ is formed by joining the feet of the medians and four  [{\it sic}] parallelograms are also formed, viz., $AEFD$, $BFED$, and $CEDF$. Let $a$, $b$, $c$; $d$, $e$, $f$ represent the three medians of $ABC$, and the three sides of $DEF$. Then the sum of the squares of the six diagonals equals the sum of the squares of the twelve sides of the parallelograms, which are equal in sets of four. That is, $a^2+b^2+c^2+d^2+e^2+f^2=4 (d^2+e^2+f^2)$ or $a^2+b^2+c^2=3 (d^2+e^2+f^2) = 3/4 (AB^2+BC^2+CA^2)$.
\end{quote}
The configuration of triangles is displayed in Figure~5 [the figure is omitted from this version]. The problem is solved using the fact that, in a parallelogram, the sum of the squares of the lengths of the two diagonals equals twice the squares of the two side lengths.

One of the most advanced topics treated in Corbin's published work is the theory of  matrices and determinants, a topic that was not usually taught at American universities in the 19th century. This topic is the subject of Corbin's only expository article,  ``Note on Elimination." This 1896  {\it Monthly} publication  \cite{corbin96e} is a short one-page note in which Corbin explains a method for solving a system of two linear equations in two variables, $x$ and $y$. The method is essentially Gaussian elimination, i.e.~ solve  for one variable by eliminating the other. He summarizes the rule by: ``The difference (sum) of the products containing $x$ ($y$) is equal to the difference (sum) of the numerical products." While this method was certainly not original to Corbin, it appears that it was not well known to American mathematicians as he presents it as an alternative to the ``determinant method" (perhaps a version of Cramer's rule). 

Corbin also published two problems on determinants.  The first appeared in a 1895 issue \cite{corbin95f} of the {\it Monthly}:\footnote{The expression we give corrects a typo in \cite{corbin95f}. In the original printed version, the bottom right-most entry in the first matrix is $s-a_n^2$, and we have changed this to $(s-a_n)^2$. To see that the entry $s-a_n^2$ is a typo, see the published solution in the April 1896 issue of the {\it Monthly}.}

\begin{quote}
Find the quotient of
\begin{gather*}
	\begin{pmatrix} 
		(s-a_1)^2 	& a_1^2 		& \dots 	& a_{1}^{2} \\
		a_{2}^{2}	&(s-a_2)^2 	& \dots 	& a_{2}^{2} \\
		a_{3}^{2}	& a_{3}^{2}	& \dots 	& a_{3}^{2}  \\
		\vdots 	&	\vdots	& \vdots	& \vdots	\\
		a_{n}^{2}	& a_{n}^{2}	& \dots 	& (s-a_n)^2
	\end{pmatrix} \\ \div 
	\begin{pmatrix} 
		s-a_1 	& a_1 		& a_1 		& \dots 	& a_{1} \\
		a_{2}	& s-a_2 		& a_2		& \dots 	& a_{2} \\
		a_{3}	& a_{3}		& s-a_3	 	& \dots 	& a_{3}  \\
		\vdots 	&	\vdots	&	\vdots	& \vdots	& \vdots	\\
		a_{n}	& a_{n}		& a_{n}	 	& \dots 	& s-a_n
	\end{pmatrix}
\end{gather*}
\end{quote}

The desired expression for the quotient is
\[
	s^{n-1} \cdot \left(s + \sum \frac{a_i^{2}}{s-2 a_i}\right) \cdot  \frac{1}{1+ \sum \frac{a_i}{s-2 a_i}}.
\]
The solution is obtained by applying matrix manipulations that preserve the determinant. First, observe that the given ratio of determinants  equals:
\begin{gather*}
	\begin{pmatrix} 
		1		&	0		& 0			& \dots	&	0	\\
		1		&	(s-a_1)^2 	& a_2^2 		& \dots 	& a_{n}^{2} \\
		1		&	a_{1}^{2}	&(s-a_2)^2 	& \dots 	& a_{n}^{2} \\
		1		&	a_{1}^{2}	& a_{2}^{2}	& \dots 	& a_{n}^{2}  \\
		\vdots	&	\vdots 	&	\vdots	& \vdots	& \vdots	\\
		1		&	a_{1}^{2}	& a_{2}^{2}	& \dots 	& (s-a_n)^2
	\end{pmatrix} \\ \div 
	\begin{pmatrix} 
		1		&	0		& 0			& 0			& \dots	&	0	\\
		1		&	s-a_1 	& a_2 		& a_3 		& \dots 	& a_{n} \\
		1		&	a_{1}	& s-a_2 		& a_3		& \dots 	& a_{n} \\
		1		&	a_{1}	& a_{2}		& s-a_3	 	& \dots 	& a_{n}  \\
		\vdots	&	\vdots 	&	\vdots	&	\vdots	& \vdots	& \vdots	\\
		1		&	a_{1}	& a_{2}		& a_{3}	 	& \dots 	& s-a_n
	\end{pmatrix}
\end{gather*}
Each column now has the property that almost all entries are the same, and thus it can be simplified by subtracting a suitable multiple of the first column. This yields:
\begin{gather*}
	\begin{pmatrix} 
		1		&	-a_{1}^2		& -a_{2}^2		& \dots	&-a_{n}^{2}	\\
		1		&	s(s-2 a_1) 	& 0 			& \dots 	& 0 \\
		1		&	0			&s (s-2 a_2) 	& \dots 	& 0  \\
		1		&	0			& 0			& \dots 	& 0  \\
		\vdots	&	\vdots 		&	\vdots	& \vdots	& \vdots	\\
		1		&	0			& 0			& \dots 	& s (s-2 a_n)
	\end{pmatrix} \\ \div 
	\begin{pmatrix} 
		1		&	-a_{1}	& -a_{2}		& -a_{3}		& \dots	&	-a_{n}	\\
		1		&	s- 2 a_1 	& 0	 		& 0			& \dots 	& 0 \\
		1		&	0		& s- 2 a_2 		& 0			& \dots 	& 0 \\
		1		&	0		& 0			& s- 2 a_3	 	& \dots 	& 0  \\
		\vdots	&	\vdots 	&	\vdots	&	\vdots	& \vdots	& \vdots	\\
		1		&	0		& 0			& 0		 	& \dots 	& s- 2 a_n
	\end{pmatrix}
\end{gather*}
Computing each determinant as the Laplace expansion along the first row, we get that the first determinant is
\begin{multline*}
s^n \prod (s-2 a_i) + a_{1}^{2} s^{n-1} \prod_{i \ne 1} (s-2 a_i) +\dots \\+ a_{j}^{2} s^{n-1} \prod_{i \ne j} (s-2 a_{i})+\dots +  a_{n}^{2} s^{n-1} \prod_{i \ne n}(s-2 a_i),
\end{multline*}
and the second is
\begin{multline*}
	\prod (s-2 a_i) +a_{1} \prod_{i \ne 1} (s-2 a_i)+ a_{2} \prod_{i \ne 2} (s-2 a_i)\\
	+\dots + a_{j} \prod_{i \ne j} (s-2a_{i})+\dots+a_{n} \prod_{i \ne n} (s-2 a_i).
\end{multline*}
Factoring out $\prod (s-2 a_{i})$, we get the desired solution.

Corbin's second published problem on determinants was his last publication. In  \cite{corbin08}, Corbin asked:
\begin{quote}
	Muir gives the following problem: 
	Prove:
	\[
		\begin{vmatrix}
			1	&	a	&	a	&	a^2	\\
			1	&	b	&	b	&	b^2	\\
			1	&	c	&	c'	& 	c c'	\\
			1	&	d	&	d'	&	d d'
		\end{vmatrix}
		=
		(a-b) \begin{vmatrix}
			1	&	a b	&	a+b \\
			1	&	c d'	&	c+d'	\\
			1	&	c' d	&	c'+d
		\end{vmatrix}
	\]
	which, of course, can be solved by finding the terms of both determinants. Is there any method of changing from one form to the other which is direct?
\end{quote}
The reference to ``Muir" is a reference to Thomas Muir's  book  {\it A Treatise on the Theory of Determinants}.  The problem appears as one of the exercises in the chapter ``Determinants in general" on the general properties of the determinant.\footnote{Specifically, it appears on page 66 of Muir's book as Exercise 27  in Exercise Set IX in Chapter II.} The solution is similar to the one for the previous problem. The first matrix is manipulated using operations that preserve the determinant until its determinant is visibly equal to the determinant of the second matrix.

\subsection*{Corbin's Legacy}
American mathematical culture underwent major changes around the time Corbin's career came to an end. Four years after his death, the Mathematical Association of America was formed, in large part to provide a permanent institution  to support the {\it American Mathematical Monthly}. The {\it Monthly} became an official journal of the MAA, and  the loosely knit readership  of math enthusiasts, which   Corbin had been a part of, began to transform into an organized group of professionals. Black mathematicians continued to be members of this community. The inaugural charter members of the MAA included James T.~Cater, a graduate of the HBCU Atlanta University and a math professor at Straight University (now part of Dillard University) \cite{monthly16}.

At the Branch Normal College, events like Corbin's removal were major setbacks to efforts to run the college as an institute of higher education. However, they did not end those efforts. Faculty continued participate in national mathematical culture. In 1938,  faculty at the college (then renamed the  Arkansas Agricultural, Mechanical \& Normal College) joined the ranks of the American Mathematical Society when college professor William Louis Fields was elected to membership \cite{AMSBull39}.\footnote{Other Branch Normal faculty who were later elected to AMS membership include Morris Edward Mosley \cite{AMSBull48} and Willie E.~Clark \cite{AMSBull50}.} Fields joined the MAA three years later  \cite[p.~60]{monthly41}.\footnote{By 1951, four additional Branch Normal faculty members (Mrs.~Willie E.~Clark, Garland~D.~Kyle, Nathan T.~Seely, Jr., and H.~B.~Young \cite{monthly50, monthly51}) had joined the MAA.}

In 1903, a year after he was  removed from the Branch Normal College, Corbin delivered an address at the annual meeting of the Arkansas Negro State Teachers' Association.  The meeting was attended by over 100 school teachers and was reported in the press \cite{newspaper03}, so it provided Corbin an  important forum for voicing his views on education. His speech reads as a rebuke of the educational philosophy represented by the decision to make Isaac Fischer principal of the Branch Normal College. Corbin said
\begin{quote}
	[A] need of the negro is a supply of men with the necessary equipment and capacity to carry on his vast national enterprises. It may be a surprising thing to many that I claim that the negro has any enterprises of this nature; but the claim is easily verified.
\end{quote}
As evidence for his claim, Corbin proceeded to give examples of achievements in professions like law, education, and farming. Were Corbin to give the speech today, one imagines that he would  also include mathematics. 

\subsection*{Added after publication}
This article was originally published as \cite{kass23}. The present version includes a more extensive bibliography.  We also provide the following additional information on three different issues.

First, in footnote 20, we observed that an obituary for Corbin states that he served as president of  ``the Ouachita Baptist College" in Camden, Arkansas, and we remarked that this appears to be an error as the college (now renamed Ouachita Baptist University) was a whites-only institute located in Arkadelphia, not Camden. It is likely that there was some confusion over the name of the institute Corbin taught at. Camden was home to a different school called the ``Ouachita Academy" or ``Ouachita Industrial Academy"  that was maintained by the Black Baptist church. This is likely the school that Corbin worked at. See  \cite{jonesBook}*{p.~135} or \cite{workBook}*{p.~159} for basic details about the academy.

Second, we stated at the beginning of the section ``Corbin's mathematical work" that Corbin's first mathematical publication appeared in a 1882 issue of the {\it Mathematical Magazine}. In fact, his name appears in a 1879 issue of the periodical {\it Educational Notes and Queries} \cite{queries79b}. There he is credited with submitting a solution to a problem \cite{queries79a} that had been posed earlier that year by Franklin Pierce Matz, a professor living in King's Mountain, North Carolina. Matz asked the following, ``Given 
\begin{align*}
	m^{2}/x^{2}+a =& m^{2}/y^{2}+b\\
	=&  \sqrt{1-m^{2} ( 1/x^{2} + 1/y^{2}) + m^{4}/x^{2} y^{2}}
\end{align*}
to find the values of $x$ and $y$ on the plan of quadratics." In the published solution, the answer is obtained by eliminating $y$ and then solving for $x$. This yields the answer
\begin{align*}
	x =& \pm m \sqrt{ (2 + a + b)/(1+b-a-a^{2}) }, \\
	y =& \pm m \sqrt{ (2 + a + b)/(1+a-b-b^{2}) }.
\end{align*}

Finally, the description of Corbin's degrees in footnote 9 is potentially misleading. We described his two A.M. degrees from Ohio University as ``honorary." The term ``honorary" is often used to describe an unearned degree given to a visitor by a university in recognition their accomplishments. (For example, the football coach W.~W.~``Woody" Hayes  was awarded an honorary  Doctor of Humanities by  The Ohio State University when he spoke at a commencement ceremony.) This  description does not apply to the degrees awarded to Corbin. Ohio University records describe both degrees as awarded ``in course." The first degree was awarded three years after Corbin received his B.A., in 1856.  This was in accordance with published degree requirements. During this period, all students who completed a B.A.~degree were then awarded an A.M.~degree after three years. It is unclear why Corbin was awarded his second  A.M. degree. He received the degree in  in 1889, more than thirty years after he had graduated from the university.

Our sources for Corbin's Ph.D.~are the OU alumni publications \cite[pp.~26--27]{ou1909} and \cite[p.~14]{ou1946}. There the degree is just described as a ``Ph.D.~[that] was conferred upon him by a Baptist institution in the South." Dr.~Gladys Turner Finney suggested  that a likely candidate for the ``Baptist Institution" is Simmons College of Kentucky. The college is affiliated with the Baptist church, and the people who founded helped found it included the husband of Corbin's sister, Henry Adams.

The author thanks Dr.~Gladys Turner Finney both for her questions about footnote 9 and for sharing her research on Corbin's degrees.



%
%
%
%
%
%
%
%
%

\bibliography{Corbin}

\end{document}